\documentclass[12pt]{article}
\usepackage{color}
\definecolor{darkblue}{rgb}{0.00,0.25,0.50}
\usepackage[colorlinks,filecolor=blue,citecolor=darkblue]{hyperref}

\usepackage{amsfonts,amssymb,amsmath,amsthm}
\usepackage{url}
\usepackage{enumerate}
\usepackage[ukrainian, russian, english]{babel}
\usepackage[cp1251]{inputenc}

\usepackage{cite}
\begin{document}
\selectlanguage{ukrainian} 
\thispagestyle{empty}

\title{}

\begin{center}
\textbf{\Large Оцінки знизу колмогоровських поперечників класів інтегралів Пуассона}
\end{center}
\vskip0.5cm
\begin{center}
А.\,С.~Сердюк, В.\,В.~Боденчук\\ \emph{\small Інститут математики НАН України, Київ}
\end{center}
\vskip0.5cm

\begin{abstract}
We expand the ranges of permissible values of $n$ ($n\in\mathbb{N}$) for which Poisson kernels $P_{q,\beta}(t)=\sum\limits_{k=1}^{\infty}q^k\cos\left(kt-\dfrac{\beta\pi}{2}\right)$, ${q\in(0,1)}$, $\beta\in\mathbb{R}$, satisfy Kushpel's condition $C_{y,2n}$. As a consequence, we obtain exact values for Kolmogorov widths in the space $C$~($L$) of classes $C_{\beta,\infty}^q$~($C_{\beta,1}^q$) of Poisson integrals generated by kernels $P_{q,\beta}(t)$ in new situations. 
It is shown that obtained here results we can't obtain by using methods of finding of exact lower bounds for widths suggested by A. Pinkus.

\vskip0.6cm

Розширено області допустимих значень параметра $n$ ($n\in\mathbb{N}$), при яких ядра Пуассона $P_{q,\beta}(t)=\sum\limits_{k=1}^{\infty}q^k\cos\left(kt-\dfrac{\beta\pi}{2}\right)$, ${q\in(0,1)}$, $\beta\in\mathbb{R}$, задовольняють умову Кушпеля $C_{y,2n}$. Як наслідок, в нових ситуаціях встановлено точні значення колмогоровських поперечників в просторі $C$~($L$) класів інтегралів Пуассона $C_{\beta,\infty}^q$~($C_{\beta,1}^q$), породжених ядрами $P_{q,\beta}(t)$. Показано, що знайдені в роботі результати неможливо одержати, використовуючи методи знаходження точних оцінок знизу поперечників, які були розвинуті А.~Пінкусом.
\end{abstract}

\vskip1cm

%%%%%%%%%%%%%%%%%%%%%%%%%%%%%%%%%%%%%%%%%%%%%%%%%%%%%%%%%%%%%%%%%%%%%%%%%

Нехай $L=L_1$ простір $2\pi$-періодичних сумовних функцій $f$ з нормою
$\|f\|_1=\int\limits_{-\pi}^{\pi}|f(t)|dt$, $L_\infty$ --- простір $2\pi$-періодичних вимірних і суттєво обмежених функцій з нормою
$\|f\|_\infty=\mathop{\rm ess\,sup}\limits_{t\in\mathbb{R}\ }|f(t)|,$
$C$ --- простір $2\pi$-періодичних неперервних функцій $f$, у якому норма задається рівністю
$\|f\|_C=\max\limits_{t\in\mathbb{R}}|f(t)|$.

Інтегралом Пуассона функції $\varphi\in L^0=\{\varphi\in L: \int\limits_{-\pi}^{\pi}\varphi(t)dt=0\}$ називають функцію $f$, що зображується у вигляді згортки
\begin{equation}\label{f}
f(x)=A+\dfrac{1}{\pi}\int\limits_{-\pi}^{\pi}P_{q,\beta}(x-t)\varphi(t)dt=A+\left(P_{q,\beta}\ast\varphi\right)(x),\; A\in\mathbb{R},
\end{equation}
де
\begin{equation}\label{P_qb}
P_{q,\beta}(t)=\sum\limits_{k=1}^{\infty}q^k\cos\left(kt-\dfrac{\beta\pi}{2}\right),\, q\in(0,1),\, \beta\in\mathbb{R}
\end{equation}
--- ядро Пуассона $P_{q,\beta}(t)$ з параметрами $q$ і $\beta$.
Функцію $\varphi\in L^0$, пов’язану із $f$ за допомогою рівності \eqref{f}, називають $(q,\beta)$-похідною функції $f$ і позначають через $f_\beta^q$ (див., наприклад, \cite[с.~302]{Stepanets_2002_1}).

Множину усіх інтегралів Пуассона вигляду \eqref{f} у випадку, коли 
\begin{equation*}
\|\varphi\|_p\leqslant1,\, p=1,\infty,
\end{equation*}
будемо позначати через $C_{\beta,p}^{q}$.

Нехай $d_m(\mathfrak{M},X)$ --- $m$-вимірний поперечник за Колмогоровим центральносиметричної множини $\mathfrak{M}$ банахового простору $X$, тобто величина, що означається рівністю
\begin{equation}\label{d_m}
d_m(\mathfrak{M},X)=\inf_{F_m\subset X}\sup_{f\in \mathfrak{M}}\inf_{u\in F_m}\|f-u\|_X,
\end{equation}
де зовнішній $\inf$ розглядається по всіх $m$-вимірних лінійних підпросторах $F_m$ із $X$.

Задача про знаходження точних значень або порядкових оцінок колмогоровських поперечників для різноманітних функціональних компактів в різноманітних функціональних просторах має багату історію, ознайомитись з якою можна по роботах \cite{Pinkus_1985,Tihomirov_1976,Kornejchuk_1987,Temlyakov, Serdyuk_1995,Serdyuk_1999,Serdyuk_2005,Kushpel_1988,Kushpel_1989, Shevaldin_1992,Nguen_1994}.

Задача про обчислення поперечників $d_m(\mathfrak{M},X)$, як правило, розпадається на дві частини. Спочатку фіксується деякий підпростір ${F_m\subset X}$, $\dim F_m=m$, і обчислюється величина
\begin{equation}\label{E_n}
E(F_m,\mathfrak{M},X)=\sup_{f\in\mathfrak{M}}\inf_{u\in F_m}\|f-u\|_X,
\end{equation}
Зрозуміло, що згідно з \eqref{d_m} і \eqref{E_n}
\begin{equation}\label{d_m_E_n}
E(F_m,\mathfrak{M},X)\geqslant d_m(\mathfrak{M},X).
\end{equation}
Потім для поперечника $d_m(\mathfrak{M},X)$ отримують оцінки знизу. Позначимо через $\mathcal{T}_{2n-1}$ підпростір тригонометричних поліномів $t_{n-1}$, порядок яких не перевищує $n-1$ і розглянемо величини найкращих наближень $E_n(C_{\beta,\infty}^{q})_C=E(\mathcal{T}_{2n-1},C_{\beta,\infty}^{q},C)$ і $E_n(C_{\beta,1}^{q})_{L}=\linebreak =E(\mathcal{T}_{2n-1},C_{\beta,1}^{q},L)$.

Як випливає з \cite{Shevaldin_1992,Krein_1938,Nikolskiy_1946,Bushanskiy_1978}, для довільних $q\in(0,1)$, $\beta\in\mathbb{R}$ і $n\in\mathbb{N}$ мають місце рівності
\begin{equation*}
E_n(C_{\beta,\infty}^{q})_C=E_n(C_{\beta,1}^{q})_{L}=\|P_{q,\beta}\ast\varphi_n\|_C=
\end{equation*}
\begin{equation}\label{E_n_rivnosti}
=\dfrac{4}{\pi}\left|\sum\limits_{\nu=0}^{\infty}\dfrac{q^{(2\nu+1)n}}{2\nu+1}\sin\left((2\nu+1)\theta_n\pi-\dfrac{\beta\pi}{2}\right)\right|,
\end{equation}
де
\begin{equation*}
\varphi_n(t)=\textnormal{sign}\sin nt,
\end{equation*}
а $\theta_n=\theta_n(q,\beta)$ --- єдиний на $[0,1)$ корінь рівняння
\begin{equation}\label{theta}
\sum\limits_{\nu=0}^{\infty}q^{(2\nu+1)n}\cos\left((2\nu+1)\theta_n\pi-\dfrac{\beta\pi}{2}\right)=0.
\end{equation}
Тому, з урахуванням \eqref{d_m_E_n}, для розв'язання задачі про точні значення вказаних поперечників залишається встановити оцінки знизу
\begin{equation}\label{dno1}
d_{2n}(C_{\beta,\infty}^q, C)\geqslant\|P_{q,\beta}\ast\varphi_n\|_C,
\end{equation}
\begin{equation}\label{dno2}
d_{2n-1}(C_{\beta,1}^q, L)\geqslant\|P_{q,\beta}\ast\varphi_n\|_C.
\end{equation}

Мета даної роботи полягає у доведенні оцінок \eqref{dno1} і \eqref{dno2} у нових, не досліджених раніше, випадках.

Вперше нерівності \eqref{dno1} і \eqref{dno2} при $q\in(0,\dfrac{1}{7}]$, $\beta\in\mathbb{Z}$ і $n\in\mathbb{N}$ встановив О.К.~Кушпель у \cite{Kushpel_1988} і \cite{Kushpel_1989}. Згодом В.Т. Шевалдін \cite{Shevaldin_1992} показав, що зазначені нерівності виконуються при $q\in(0, q(\beta)]$, $\beta\in\mathbb{R}$
і $n\in\mathbb{N}$, де $q(\beta)=0{,}2$ при $\beta\in\mathbb{Z}$ і $q(\beta)=0{,}196881$ при $\beta\in\mathbb{R\setminus Z}$. 
З роботи Нгуен Тхи Тх’єу Хоа \cite[с.~211]{Nguen_1994} випливає справедливість \eqref{dno1} і \eqref{dno2} для будь-яких $q\in(0,1)$,
 $\beta=2kl$, $l\in\mathbb{Z}$, при натуральних $n$, більших деякого, залежного від $q$, номера $n_*$
(при цьому було доведено існування номера $n_*$, а конструктивного способу знаходження $n_*$ по $q$ не вказано). 
Авторами \cite{My_arxiv_2012} доведено нерівності \eqref{dno1} і \eqref{dno2} для довільних $q\in(0,1)$, $\beta\in\mathbb{R}$
і $n\in\mathbb{N}$, $n\geqslant n_q$, де $n_q$ --- найменший з номерів $n\geqslant9$  при фіксованому $q\in(0,1)$, для яких виконується нерівність
\begin{equation*}
\dfrac{43}{10(1-q)}q^{\sqrt{n}}+\dfrac{160}{57(n-\sqrt{n})}\; \dfrac{q}{(1-q)^2}\leqslant
\end{equation*}
\begin{equation}\label{umova_n_0_1}
\leqslant
\left(\dfrac{1}{2}+\dfrac{2q}{(1+q^2)(1-q)}\right)\left(\dfrac{1-q}{1+q}\right)^{\frac {4}{1-q^2}}.
\end{equation}

У даній роботі область допустимих значень параметра $n$, для яких справджуються оцінки \eqref{dno1} і \eqref{dno2} вдалось дещо розширити. Перш ніж сформулювати основний результат роботи, розглянемо при кожному фіксованому $q\in(0,1)$ найменший з номерів $n\geqslant9$, для яких виконується нерівність
\begin{equation*}
\dfrac{43}{10(1-q)}q^{\sqrt{n}}+\dfrac{q}{(1-q)^2} \min\left\{\dfrac{160}{57(n-\sqrt{n})}, \dfrac{8}{3n-7\sqrt{n}}\right\}
 \leqslant
\end{equation*}
\begin{equation}\label{umova_n_00}
\leqslant
\left(\dfrac{1}{2}+\dfrac{2q}{(1+q^2)(1-q)}\right)\left(\dfrac{1-q}{1+q}\right)^{\frac {4}{1-q^2}}.
\end{equation}
Будемо позначати цей номер через $n_q^{*}$.

\textbf{Теорема 1.}
\emph{Нехай $q\in(0,1)$.
Тоді для довільного $\beta\in \mathbb{R}$ та всіх номерів $n\geqslant n_q^{*}$ мають місце оцінки \eqref{dno1} та \eqref{dno2}.
}

Зауважимо, що згідно з \eqref{umova_n_0_1} і \eqref{umova_n_00} $n_q^{*}\leqslant n_q$, тому з теореми~1 випливає теорема~2 роботи \cite{My_arxiv_2012}.
Отже, сформульована теорема у порівнянні з теоремою~2 роботи \cite{My_arxiv_2012} дозволяє розширити область допустимих значень параметра $n$, для яких виконуються оцінки \eqref{dno1} і \eqref{dno2}, при тих значеннях $q$, для яких $n_q>n_q^{*}$. Обчислення показують, що нерівність $n_q>n_q^{*}$ виконується зокрема при усіх $q\in[0{,}4925,\,1)$.
Наприклад, при $q=0{,}5$ $n_q=969$, а $n_q^*=963$.

\textbf{\emph{Доведення.}} Будемо використовувати запропонований О.К.~Кушпелем метод знаходження оцінок знизу колмогоровських поперечників класів згорток із твірними ядрами, що задовольняють так звану умову $C_{y, 2n}$.
Наведемо необхідні означення. 
Нехай $\Delta_{2n}=\{0=x_0<x_1<\dots<x_{2n}=2\pi\}$, $x_k=k\pi/n$ --- розбиття проміжку $[0,2\pi]$ і
\begin{equation}\label{Ps}
P_{q,\beta,1}(t)=(P_{q,\beta}\ast B_{1})(t)=\sum\limits_{k=1}^{\infty}\frac{q^k}{k}\cos\left(kt-\dfrac{(\beta+1)\pi}{2}\right),
\end{equation}
\begin{equation*}
q\in(0,1),\; \beta\in\mathbb{R},
\end{equation*}
де $B_{1}(t)=\sum\limits_{k=1}^{\infty}k^{-1}\sin kt$ --- ядро Бернуллі.
Фундаментальним $SK$-сплайном називають функцію $\overline{SP}_{q,\beta,1}(\cdot)=\overline{SP}_{q,\beta,1}(y, \cdot)$ виду
\begin{equation}\label{SK}
SP_{q,\beta,1}(\cdot)=\alpha_{0}+\sum\limits_{k=1}^{2n} \alpha_k P_{q,\beta,1}(\cdot-x_k),\;\sum\limits_{k=1}^{2n} \alpha_k=0,
\end{equation}
\begin{equation*}
\alpha_k\in\mathbb{R},\; k= 0, 1, \dots, 2n,
\end{equation*}
що задовольняє співвідношення
\begin{equation*}
\overline{SP}_{q,\beta,1}(y, y_k)=\delta_{0,k}=
\begin{cases}
0, & k=\overline{1,2n-1},    \\
1, & k=0,
\end{cases}
\end{equation*}
де $y_k=x_k+y$, $x_k=k\pi/n$, $y\in[0,\dfrac{\pi}{n})$. 

В силу \eqref{Ps} і згідно з означенням поняття $(q,\beta)$-похідної
\begin{equation}\label{e15}
(P_{q,\beta,1}(\cdot))_\beta^q=B_{1}(\cdot),
\end{equation}
тому з \eqref{SK} маємо
\begin{equation}\label{e16}
(SP_{q,\beta,1}(\cdot))_\beta^q=\sum\limits_{k=1}^{2n}\alpha_kB_{1}(\cdot-x_k),\; \sum_{k=1}^{2n}\alpha_k=0.
\end{equation}
Рівності в \eqref{e15} і \eqref{e16} слід розуміти як рівності функцій з $L_1$, тобто майже скрізь. В силу леми 2.3.4 роботи \cite[с. 76]{Kornejchuk_1987} функція, що знаходиться в правій частині рівності \eqref{e16} є константою на кожному інтервалі $(x_k, x_{k+1})$.
Тому, серед $(q, \beta)$-похідних будь-якого сплайна вигляду \eqref{SK}, а значить, і для фундаментального сплайна $\overline{SP}_{q,\beta,1}(\cdot)$, існує функція, яка є сталою на кожному інтервалі $(x_k, x_{k+1})$. Надалі саме таку функцію будемо розуміти під записом $(\overline{SP}_{q,\beta,1}(\cdot))_\beta^q$.

\textbf{Означення 1.} \emph{Будемо казати, що для деякого дійсного числа $y$ і розбиття $\Delta_{2n}$ ядро $P_{q,\beta}(\cdot)$ вигляду \eqref{P_qb} задовольняє умову $C_{y,2n}$ (і записувати $P_{q,\beta}\in C_{y,2n}$), якщо для цього ядра існує єдиний фундаментальний сплайн $\overline{SP}_{q,\beta,1}(y,\cdot)$ і для нього виконуються рівності
\begin{equation*}
\emph{\textrm{sign}}(\overline{SP}_{q,\beta,1}(y,t_k))_\beta^q=(-1)^k\varepsilon e_k,\, k=\overline{0,2n-1},
\end{equation*}
де $t_k=(x_{k}+x_{k+1})/2,$ $e_k$ дорівнює або 0, або 1, а $\varepsilon$ приймає значення $\pm1$ і не залежить від $k$.}

Із робіт О.К.~Кушпеля \cite{Kushpel_1988,Kushpel_1989} випливає, що якщо ядро Пуассона $P_{q,\beta}$ задовольняє умову $C_{y_0,2n}$, де $y_0$ --- точка, в якій функція ${|P_{q,\beta}\ast\varphi_n|}$, досягає найбільшого значення, тобто $|(P_{q,\beta}\ast\varphi_n)(y_0)|=\linebreak=\|P_{q,\beta}\ast\varphi_n\|_C$, то для поперечників класів згорток з таким ядром мають місце оцінки \eqref{dno1} та \eqref{dno2}.
З урахуванням того, що ${y_0=y_0(n,q,\beta)=\dfrac{\theta_n\pi}{n}}$, де $\theta_n$ --- корінь рівняння \eqref{theta}, $\theta_n\in[0,1)$, то для доведення теореми~1 достатньо встановити включення $P_{q,\beta}\in C_{y_0, 2n}$. 

Відповідно до леми~1 роботи \cite{My_arxiv_2012} для довільного ${t\in(\dfrac{(k-1)\pi}{n},\dfrac{k\pi}{n})}$, $k=\overline{1,2n}$, має місце представлення
\begin{equation}\label{SP_Phi_}
(\overline{SP}_{q,\beta,1}(y_0,t))_\beta^q=
(-1)^{k+s+1}\frac{\pi}{4nq^{n}}\;
(\mathcal{P}_q(t_k-y_0)+\sum_{m=1}^5\gamma_m(y_0)),
\end{equation}
в якій $\mathcal{P}_q(t)$ --- ядро Пуассона рівняння теплопровідності
\begin{equation*}
\mathcal{P}_q(t)= \dfrac{1}{2}+2\sum_{j=1}^{\infty}\dfrac{\cos jt}{q^j+q^{-j}},
\end{equation*}
\textit{а}
\begin{align}
\label{gamma_1}
\gamma_1(y_0)&=\gamma_1(k,y_0)=2\sum_{j=[\sqrt{n}]+1}^{n-1}\dfrac{\cos j(t_k-y_0)}{\dfrac{n}{q^{n}}|\lambda_{n-j}(y_0)|\cos\dfrac{j\pi}{2n}},
\\ 
\gamma_2(y_0)&=\gamma_2(k,y_0)= \nonumber
\\
\label{gamma_2}
&=(-1)^{s}\frac{q^{n}}{n}
\left(\dfrac{z_{0}(y_0)}{|\lambda_{n}(y_0)|^2}+
2\sum_{j=1}^{n-1}\dfrac{z_{j}(y_0)}{|\lambda_{n-j}(y_0)|^2\cos\dfrac{j\pi}{2n}}\right),
\\
\label{gamma_3}
\gamma_3(y_0)&=-\dfrac{R_0(y_0)\dfrac{n}{q^{n}}}{2(2+R_0(y_0)\dfrac{n}{q^{n}})},
\\
\label{gamma_4}
\gamma_4(y_0)&=\gamma_4(k,y_0)=-2\sum_{j=1}^{[\sqrt{n}]}\dfrac{\delta_{j}(y_0)\cos j(t_k-y_0)}{\dfrac{n}{q^{n}}|\lambda_{n-j}(y_0)|\cos\dfrac{j\pi}{2n}},
\\
\label{gamma_5}
\gamma_5(y_0)&=\gamma_5(k,y_0)=-2\sum\limits_{j=[\sqrt{n}]+1}^{\infty}\dfrac{\cos j(t_k-y_0)}{q^j+q^{-j}},
\\
\label{delta_0}
\delta_{j}(y_0)&=\dfrac{n|\lambda_{n-j}(y_0)|\cos\dfrac{j\pi}{2n}}{(q^{-j}+q^{j})q^n}-1,\;j=\overline{0,[\sqrt{n}]},
\\
z_{j}(y_0)&=|r_{j}(y_0)|\cos(j(t_k-y_0)+\arg(r_{j}(y_0)))+\nonumber
\\
\nonumber
&+(-1)^{s+1}R_{j}(y_0)\cos(j(t_k-y_0)),\;j=\overline{0,n-1},
\\
\nonumber
R_{j}(y_0)&=|\lambda_{n-j}(y_0)|- \dfrac{q^{n-j}}{n-j}-\dfrac{q^{n+j}}{n+j},\;j=\overline{0,n-1},
\end{align}
\begin{equation*}
\lambda_{n-j}(y_0)=e^{-ijy_0}\left((-1)^s\left(\dfrac{q^{n-j}}{n-j}+\dfrac{q^{n+j}}{n+j}\right)+r_{j}(y_0)\right),
j=\overline{0,n-1},
\end{equation*}
\begin{align}
\nonumber
&r_{j}(y_0)=\sum_{\nu=1}^3 r_{j}^{(\nu)}(y_0),
\\
\nonumber
&r_{j}^{(1)}(y_0)=
\dfrac{q^{3n-j}e^{i(3ny_0-\frac{(\beta+1)\pi}{2})}}{3n-j}+
\\
\nonumber
&\hspace{50px}+\sum\limits_{m=2}^{\infty}\left(\dfrac{q^{(2m+1)n-j}e^{i((2m+1)ny_0-\frac{(\beta+1)\pi}{2})}}{(2m+1)n-j}+\right.
\\
\label{r_nj_1}
\nonumber
&\hspace{50px}\left.+\dfrac{q^{(2m-1)n+j}e^{-i((2m-1)ny_0-\frac{(\beta+1)\pi}{2})}}{(2m-1)n+j}\right),
\\
\nonumber
&r_{j}^{(2)}(y_0)=i\left(\dfrac{q^{n+j}}{n+j}-\dfrac{q^{n-j}}{n-j}\right)\cos(ny_0-\frac{\beta\pi}{2}),
\\
\nonumber
&r_{j}^{(3)}(y_0)=(-1)^s\left(\dfrac{q^{n-j}}{n-j}+\dfrac{q^{n+j}}{n+j}\right)(|\sin(ny_0-\dfrac{\beta\pi}{2})|-1)
\end{align}
$t_k=\dfrac{k \pi}{n}-\dfrac{\pi}{2n}$, а величина $s=s(n,q,\beta)$ означена рівністю
\begin{equation*}
(-1)^s=\mathop{\mathrm{sign}} \sin(ny_0-\dfrac{\beta \pi}{2}).
\end{equation*}

Згідно з лемою~2 роботи \cite{My_arxiv_2012} для довільного $x\in\mathbb{R}$ і довільного $q\in(0,1)$
\begin{equation}\label{f_x}
\mathcal{P}_q(x)>\left(\dfrac{1}{2}+\dfrac{2q}{(1+q^2)(1-q)}\right)\left(\dfrac{1-q}{1+q}\right)^{\frac {4}{1-q^2}}.
\end{equation}

В силу формул \eqref{SP_Phi_} і \eqref{f_x} включення $P_{q,\beta}\in C_{y_0, 2n}$ буде доведене для всіх $\beta\in\mathbb{R}$,  якщо вдасться встановити справедливість наступної нерівності:
\begin{equation}\label{gs0}
\left(\dfrac{1}{2}+\dfrac{2q}{(1+q^2)(1-q)}\right)\left(\dfrac{1-q}{1+q}\right)^{\frac {4}{1-q^2}}+\sum\limits_{k=1}^5\gamma_k(y_0)\geqslant0,
\end{equation}
де величини $\gamma_k(y_0)$, $k=\overline{1,5}$, задані рівностями \eqref{gamma_1}--\eqref{gamma_5}. 
Оцінку зверху суми $\sum\limits_{k=1}^5|\gamma_k(y_0)|$ дає наступна лема, що є узагальненням леми~3 роботи \cite{My_arxiv_2012}.

\textbf{Лема 1.} \textit{Нехай $q\in(0,1)$, $\beta\in\mathbb{R}$, ${y_0=y_0(n,q,\beta)=\dfrac{\theta_n\pi}{n}}$, де $\theta_n$ --- корінь рівняння \eqref{theta} і $\theta_n\in[0,1)$, а величини $\gamma_k(y_0)$, $k=\overline{1,5}$ задаються рівностями \eqref{gamma_1}--\eqref{gamma_5}. 
Тоді при $n\geqslant9$ та при виконанні умови}
\begin{equation}\label{umova_z}
\dfrac{q^{n}}{1-q^{2n}}\leqslant \dfrac {7q^{\sqrt{n}}}{37n^2}
\end{equation}
\textit{справедлива оцінка}
\begin{equation}\label{g}
\sum\limits_{k=1}^5|\gamma_k(y_0)|\leqslant
\dfrac{43}{10(1-q)}q^{\sqrt{n}}+\dfrac{q}{(1-q)^2} \min\left\{\dfrac{160}{57(n-\sqrt{n})}, \dfrac{8}{3n-7\sqrt{n}}\right\}.
\end{equation}

\textbf{Доведення.} В ході доведення леми~3 роботи \cite{My_arxiv_2012} для $n\geqslant9$ було встановлено, що
\begin{equation}\label{sum_gi}
|\gamma_1(y_0)|+|\gamma_2(y_0)|+|\gamma_3(y_0)|+|\gamma_5(y_0)|\leqslant \dfrac{43}{10(1-q)}q^{\sqrt{n}},
\end{equation}
а також, що 
 \begin{equation}\label{g4_0}
 |\gamma_4(y_0)|
<\dfrac{160}{57(n-\sqrt{n})}\; \dfrac{q}{(1-q)^2}.
 \end{equation}
 
Отже, щоб одержати \eqref{g}, досить показати, що при $n\geqslant9$ для величини $\gamma_4(y_0)$ виконується також нерівність 
\begin{equation}\label{g4}
 |\gamma_4(y_0)|<\dfrac{8}{3n-7\sqrt{n}}\; \dfrac{q}{(1-q)^2}.
 \end{equation}
 Згідно з формулою (83) роботи \cite{My_arxiv_2012}
\begin{equation}\label{delta}
|\delta_{j}(y_0)|\leqslant \dfrac{4j}{3(n-j)}.
\end{equation}
Записавши рівність \eqref{delta_0} у вигляді
\begin{equation*}
\dfrac{n}{q^{n}}|\lambda_{n-j}(y_0)|\cos\dfrac{j\pi}{2n}=(q^j+q^{-j})(1+\delta_{j}(y_0)),
\end{equation*}
з \eqref{gamma_4} та \eqref{delta} одержуємо, що при $n\geqslant9$
 \begin{equation*}
 |\gamma_4(y_0)|\leqslant
2\sum_{j=1}^{[\sqrt{n}]}\dfrac{\dfrac{4j}{3(n-j)}}{|1-\dfrac{4j}{3(n-j)}|}q^j=2\sum_{j=1}^{[\sqrt{n}]}\dfrac{4j}{3n-7j}q^j\leqslant
 \end{equation*}
 \begin{equation*}
 \leqslant
\dfrac{8}{3n-7\sqrt{n}}\sum_{j=1}^{[\sqrt{n}]}jq^j<\dfrac{8}{3n-7\sqrt{n}}\sum_{j=1}^{\infty}jq^j<
 \end{equation*}
 \begin{equation*}
<\dfrac{8}{3n-7\sqrt{n}}\; \dfrac{q}{(1-q)^2}.
 \end{equation*}
Тим самим \eqref{g4} доведено. З \eqref{sum_gi}, \eqref{g4_0} та \eqref{g4} одержуємо \eqref{g}. Лему доведено.

З леми~1 випливає, що при $n\geqslant9$ за умов \eqref{umova_n_00} та \eqref{umova_z} справедлива нерівність \eqref{gs0}, а отже, як наслідок, і оцінки \eqref{dno1} та \eqref{dno2}.
Тому для остаточного доведення теореми залишилось показати, що при $n\geqslant9$ нерівність \eqref{umova_n_00} забезпечує виконання умови \eqref{umova_z}.

Зазначимо, що при $q\in(0,\dfrac{91}{250}]$ умова \eqref{umova_z} виконується для довільних ${n\geqslant9}$.
Для того, щоб у цьому переконатись досить помітити, що послідовність $\xi(n)=(n-\sqrt{n})\ln\dfrac{91}{250}+2\ln n-\ln\left(\dfrac{7}{37}\left(1-\left(\dfrac{91}{250}\right)^{18}\right)\right)$ монотонно спадна при $n\geqslant9$ і $\xi(9)<0$. Тому при $n\geqslant9$
\begin{equation}\label{q003}
(n-\sqrt{n})\ln\dfrac{91}{250}+2\ln n-\ln\left(\dfrac{7}{37}\left(1-\left(\dfrac{91}{250}\right)^{18}\right)\right)<0.
\end{equation}
Нерівність \eqref{q003} еквівалентна нерівності
\begin{equation*}
\frac{\left(\dfrac{91}{250}\right)^{n-\sqrt{n}}}{1-\left(\dfrac{91}{250}\right)^{18}}<\frac{7}{37n^2},
\end{equation*}
а, тому при $q\in(0,\dfrac{91}{250}]$
\begin{equation*}
\frac{q^{n-\sqrt{n}}}{1-q^{2n}}<\frac{\left(\dfrac{91}{250}\right)^{n-\sqrt{n}}}{1-\left(\dfrac{91}{250}\right)^{18}}<\frac{7}{37n^2}.
\end{equation*}

Отже, для доведення теореми~1 достатньо показати, що при $n\geqslant9$ і ${q\in(\dfrac{91}{250},1)}$  має місце імплікація
\begin{equation}\label{imp1}
\eqref{umova_n_00}\Rightarrow\eqref{umova_z}. 
\end{equation} 

Для номерів $n$ таких, що
\begin{equation*}
\min\left\{\dfrac{8}{3n-7\sqrt{n}}, \dfrac{160}{57(n-\sqrt{n})}\right\}=\dfrac{160}{57(n-\sqrt{n})}
\end{equation*}
імплікація \eqref{imp1} доведена у \cite{My_arxiv_2012}. Тому залишається довести її при тих $n\in\mathbb{N}$, для яких
\begin{equation}\label{min}
\min\left\{\dfrac{8}{3n-7\sqrt{n}}, \dfrac{160}{57(n-\sqrt{n})}\right\}=\dfrac{8}{3n-7\sqrt{n}}.
\end{equation}

Оскільки 
\begin{equation*}
\dfrac{1}{2}+\dfrac{2q}{(1+q^2)(1-q)}< \dfrac{1+q}{1-q},\; q\in(0,1),
\end{equation*}
то з \eqref{umova_n_00} та \eqref{min} випливає нерівність
 \begin{equation*}
\dfrac{8}{3n-7\sqrt{n}}\; \dfrac{q}{(1-q)^2}<\left(\dfrac{1-q}{1+q}\right)^{\frac {4}{1-q^2}-1},
\end{equation*}
а, отже, й еквівалентна їй нерівність
\begin{equation}\label{sqrt_n0}
3n-7\sqrt{n}-\dfrac{8 q}{(1-q)^2} \left(\dfrac{1+q}{1-q}\right)^{\frac {4}{1-q^2}-1}>0.
\end{equation}
З \eqref{sqrt_n0} випливає
\begin{equation}\label{n1}
n>\dfrac{8q}{3(1-q)^2} \left(\dfrac{1+q}{1-q}\right)^{3}.
\end{equation}

Отже, при $n\geqslant9$ і ${q\in(0,1)}$
\begin{equation}\label{imp2}
\eqref{umova_n_00}\Rightarrow\eqref{n1}.
\end{equation}

Далі покажемо, що при $n\geqslant9$ і ${q\in(0,1)}$ нерівність \eqref{umova_z} випливає з нерівності
\begin{equation}\label{n2}
n>\left(\frac{9(1+q)}{4(1-q)}\right)^{2}.
\end{equation}
Оскільки (див., наприклад, \cite[с. 58]{Gradshteyn_1963}) для довільного $q\in(0,1)$
\begin{equation*}
\ln \frac{1}{q}=2\sum\limits_{k=1}^\infty \frac{1}{2k-1}\left(\frac{1-q}{1+q}\right)^{2k-1}>2\frac{1-q}{1+q},
\end{equation*}
то
\begin{equation}\label{eq:104}
\left(\frac{9(1+q)}{4(1-q)}\right)^{2}>\left(\frac{9}{4\frac{1-q}{1+q}}\right)^{\frac {125}{79}}>\left(\frac{9}{2\ln 1/q}\right)^{\frac {125}{79}}.
\end{equation}
Із \eqref{n2} і \eqref{eq:104} випливає нерівність
\begin{equation*}
n>\left(\frac{9}{2\ln 1/q}\right)^{\frac {125}{79}},
\end{equation*}
яка еквівалентна нерівності
\begin{equation}\label{eq:99}
\frac 23 n\ln \frac 1q>3 n^{\frac{46}{125}}.
\end{equation}

Оскільки при $n\in\mathbb{N}\;$ $\ln n<n^{\frac {46}{125}}$ і при $n\geqslant9\;$ $1-\dfrac{1}{\sqrt{n}}\geqslant\dfrac{2}{3}$, то з \eqref{eq:99} випливає
\begin{equation}\label{eq:101}
n\left(1-\frac{1}{\sqrt{n}}\right)\ln \frac 1q>3\ln n.
\end{equation}

При $n\geqslant9$ із \eqref{eq:101} одержуємо
\begin{equation*}
\dfrac{1}{q^{n}}> \dfrac{n^3}{q^{\sqrt{n}}}>\dfrac{9n^2}{q^{\sqrt{n}}}>\dfrac{38n^2}{7q^{\sqrt{n}}}=\dfrac{37n^2}{7q^{\sqrt{n}}}+\dfrac{n^2}{7q^{\sqrt{n}}}>\dfrac{37n^2}{7q^{\sqrt{n}}}+q^{n}.
\end{equation*}

Отже, при $n\geqslant9$ і $q\in(0, 1)$
\begin{equation}\label{imp3}
\eqref{n2}\Rightarrow\eqref{umova_z}.
\end{equation}

Залишилось довести, що при ${q\in(\dfrac{91}{250},1)}$ і $n\geqslant9$ 
\begin{equation}\label{imp4}
\eqref{n1}\Rightarrow\eqref{n2}.
\end{equation}
Для цього розглянемо різницю $v(q)$ правих частин в нерівностях \eqref{n1} та \eqref{n2}
\begin{equation*}
v(q)=\dfrac{8q}{3(1-q)^2} \left(\dfrac{1+q}{1-q}\right)^{3}-\left(\frac{9(1+q)}{4(1-q)}\right)^{2}=
\end{equation*}
\begin{equation}\label{yq}
=
\left(\dfrac{1+q}{1-q}\right)^{2}\left(\dfrac{8q(1+q)}{3(1-q)^3} -\left(\frac{9}{4}\right)^{2}\right).
\end{equation}

Оскільки ${q\in(\dfrac{91}{250},1)}$, то
\begin{equation}\label{yq1}
\dfrac{8q(1+q)}{3(1-q)^3} -\left(\frac{9}{4}\right)^{2} >0.
\end{equation}
З \eqref{yq} та \eqref{yq1} отримуємо нерівність $v(q)>0$, а разом з нею і \eqref{imp4}.
Об’єднуючи формули \eqref{imp2}, \eqref{imp3} та \eqref{imp4} одержуємо \eqref{imp1} при ${q\in(\dfrac{91}{250},1)}$.
Теорему доведено.

Зрозуміло, що при $q\in(q(\beta),1)$ знайдені у теоремі~1 оцінки знизу для колмогоровських поперечників не випливають з відомих раніше результатів, отриманих в \cite{Kushpel_1988,Shevaldin_1992,Nguen_1994}. 
Покажемо на прикладі ядер Пуассона $P_{q,0}(t)$ та $P_{q,1}(t)$ при $q=0{,}21$, 
що їх також неможливо отримати, користуючись методами і підходами, які розвинуто А.~Пінкусом \cite{Pinkus_1985} для класів згорток із ядрами, які не збільшують осцилляції.
З цією метою наведемо деякі означення та твердження. Задамо не нульовий вектор $\mathbf{x}=(x_1, \dots, x_n)$, $x_i\in\mathbb{R}$. 
Позначимо через $S(\mathbf{x})$ число змін знаку в послідовності $x_1, \dots, x_n$ без урахування нульових членів, а через $S_c(\mathbf{x})$ --- число циклічних змін знаку в $\mathbf{x}$, тобто
\begin{equation*}
S_c(\mathbf{x})=\max_i S(x_i, x_{i+1}, \dots, x_n, x_1, \dots, x_i)=
\end{equation*}
\begin{equation*}
=S(x_k, x_{k+1}, \dots, x_n, x_1, \dots, x_k),
\end{equation*}
де $k$ --- довільне ціле число для якого $x_k\not=0$. 
Для кусково-неперервної дійснозначної $2\pi$-періодичної функції $f(x)$ позначимо ${S_c(f)=\sup S_c(f(x_1), \dots, f(x_m))}$, де $m\in\mathbb{N}$, а супремум розглядається по всіх $x_1<\dots<x_m<x_1+2\pi$.

\textbf{Означення 2.} \emph{Неперервну дійснозначну $2\pi$-періодичну функцію $K(\cdot)$ називають $\emph{\text{CVD}}_{2n}$-ядром (ядром, що не збільшує осциляції) і позначають $K\in\emph{\text{CVD}}_{2n}$, якщо виконується нерівність
\begin{equation*}
S_c(K\ast f)\leqslant S_c(f),
\end{equation*}
для всіх $f$ таких, що ${S_c(f)\leqslant 2n}$.}

Кажуть, що ядро $\phi(x)$ є циклічним ядром частот Поліа порядку $2n+1$ і позначають $\phi\in\text{CPF}_{2n+1}$, якщо 
\begin{equation*}
D_{2l+1}(\mathbf{x},\mathbf{y})=\det(\phi(x_i-y_j))_{i,j=1}^{2l+1}\geqslant0,
\end{equation*}
де $0\leqslant x_1<\dots<x_{2l+1}<2\pi$, $0\leqslant y_1<\dots<y_{2l+1}<2\pi$, $l=0, 1, \dots, n$.

Співвідношення між $\text{CPF}_{2n+1}$ та $\text{CVD}_{2n}$ ядрами містяться у наступному твердженні, що належить Мерхюберу, Шонбергу та Вільямсону \cite{Mairhuber} (див. також \cite[с.~67]{Pinkus_1985}).

\textbf{Лема 2.} \emph{Нехай $\phi(x)\in C$ та $\phi(x)$ має ранг не менший за $2n+2$, тобто існує розбиття $y_i$, $i=\overline{1,2n+2}$, проміжка $[0, 2\pi)$ таке, що $0\leqslant y_1<\dots<y_{2n+2}<2\pi$ і для якого $\dim(\emph{\text{span}}\{\phi(x-y_i)\}_{i=1}^{2n+2})=2n+2$. Тоді $\phi(x)\in\emph{\text{CVD}}_{2n}$ тоді і тільки тоді, коли $\varepsilon\phi(x)\in\emph{\text{CPF}}_{2n+1}$ для деякого фіксованого $\varepsilon=\pm1$.}

Як випливає із леми~1.3 роботи \cite{Kushpel_1985}, система функцій ${\{P_{q,\beta}(x-y_i)\}_{i=1}^{2n+2}}$ лінійно незалежна і, отже, ${\dim(\text{span}\{P_{q,\beta}(x-y_i)\}_{i=1}^{2n+2})=2n+2}$. 
Тому згідно з лемою 2, щоб довести, що ядра Пуассона $P_{q_0,\beta_k}(t)$ при $q_0=0{,}21$ і $\beta_1=0$ та $\beta_2=1$ не є $\text{CVD}_{2n}$-ядрами ні при яких $n\in\mathbb{N}$, достатньо показати, що знайдуться вектори  $\mathbf{x}=(x_{1}, x_{2}, x_{3})$, $0\leqslant x_{1}<x_{2}<x_{3}<2\pi$, та $\mathbf{y}=(y_{1}, y_{2}, y_{3})$, $0\leqslant y_{1}<y_{2}<y_{3}<2\pi$, для яких детермінант $D_{3}(\mathbf{x},\mathbf{y})$ змінює знак. Виберемо вектори 
 $\mathbf{x}^{(k)}=(x_{1}^{(k)}, x_{2}^{(k)}, x_{3}^{(k)})$ та $\mathbf{y}^{(k)}=(y_{1}^{(k)}, y_{2}^{(k)}, y_{3}^{(k)})$, $k=1, 2$, наступним чином:
\begin{equation*}
x_{1}^{(1)}=\frac{\pi}{18},\; x_{2}^{(1)}=\frac{\pi}{9},\; x_{3}^{(1)}=\frac{\pi}{6},\; 
y_{1}^{(1)}=\frac{13\pi}{36},\; y_{2}^{(1)}=\frac{11\pi}{30},\; y_{3}^{(1)}=\frac{67\pi}{180},
\end{equation*}
\begin{equation*}
x_{1}^{(2)}=\frac{\pi}{18},\; x_{2}^{(2)}=\frac{\pi}{9},\; x_{3}^{(2)}=\frac{\pi}{6},\; 
y_{1}^{(2)}=\frac{13\pi}{30},\; y_{2}^{(2)}=\frac{10\pi}{9},\; y_{3}^{(2)}=\frac{7\pi}{6}.
\end{equation*}
Обчислення показують, що для ядра $P_{q_0,0}$
\begin{equation*}
D_{3}(\mathbf{x}^{(1)},\mathbf{y}^{(1)})< -9{,}98\cdot 10^{-10},\;
D_{3}(\mathbf{x}^{(2)},\mathbf{y}^{(2)})> 1{,}97\cdot 10^{-6},
\end{equation*}
а для ядра $P_{q_0,1}$ 
\begin{equation*}
D_{3}(\mathbf{x}^{(1)},\mathbf{y}^{(1)})< -1{,}3\cdot 10^{-8},\;
D_{3}(\mathbf{x}^{(2)},\mathbf{y}^{(2)})> 1{,}17\cdot 10^{-6}.
\end{equation*}
Отже, в силу леми~2 для будь-яких $n\in\mathbb{N}$ $P_{q_0,\beta_k}(t)\not\in\text{CVD}_{2n}$, $\beta_1=0$, $\beta_2=1$. 

Теорема~1, формули \eqref{d_m_E_n} і \eqref{E_n_rivnosti}, а також теорема~2 роботи В.Т.~Шевалдіна \cite{Shevaldin_1992} дають змогу записати наступне твердження про точні значення поперечників $d_{m}(C_{\beta,\infty}^q,C)$ та $d_{2m-1}(C_{\beta,1}^q,L)$, яке охоплює відомі на даний час результати \cite{Kushpel_1988,Shevaldin_1992,Nguen_1994,My_arxiv_2012,My_Dop_2013}. Для його формулювання позначимо
\begin{equation*}
n_{q,\beta}=
\begin{cases}
1,&\hspace{-7pt}\text{якщо $q\in(0, 0{,}2]$ і $\beta\in\mathbb{Z}$ або $q\in(0,0{,}196881]$ і $\beta\in\mathbb{R\setminus Z}$,}
\\
n_q^{*},&\hspace{-7pt}\text{якщо $q\in(0{,}2,1)$ і $\beta\in\mathbb{Z}$ або $q\in(0{,}196881, 1)$ і $\beta\in\mathbb{R\setminus Z}$.}
\end{cases}
\end{equation*}

\textbf{Теорема 2.}
\emph{Нехай $q\in(0,1)$.
Тоді для довільного $\beta\in \mathbb{R}$ та усіх номерів $n\geqslant n_{q,\beta}$
мають місце рівності
\begin{equation*}
d_{2n}(C_{\beta,\infty}^{q},C)=d_{2n-1}(C_{\beta,\infty}^q,C)= d_{2n-1}(C_{\beta,1}^q,L)
=
\end{equation*}
\begin{equation*}
=E_n(C_{\beta,\infty}^{q})_C
=E_n(C_{\beta,1}^{q})_L
=\|P_{q,\beta}\ast\varphi_n\|_C
=
\end{equation*}
\begin{equation}\label{dn}
=\dfrac{4}{\pi}\left|\sum\limits_{\nu=0}^{\infty}\dfrac{q^{(2\nu+1)n}}{2\nu+1}\sin\left((2\nu+1)\theta_n\pi-\dfrac{\beta\pi}{2}\right)\right|,
\end{equation}
де $\theta_n=\theta_n(q,\beta)$ --- єдиний на $[0,1)$ корінь рівняння \eqref{theta}. }

Теорема~2 дозволяє оцінити асимптотичну при ${n\to\infty}$ поведінку поперечників $d_{2n}(C_{\beta,\infty}^{q},C)$, $d_{2n-1}(C_{\beta,\infty}^q,C)$ та $d_{2n-1}(C_{\beta,1}^q,L)$.

\textbf{Теорема 3.}
\emph{Нехай $q\in(0,1)$ та $\beta\in \mathbb{R}$.
Тоді при $n\geqslant n_{q,\beta}$ }
\begin{equation*}
d_{2n}(C_{\beta,\infty}^{q},C)=d_{2n-1}(C_{\beta,\infty}^q,C)= d_{2n-1}(C_{\beta,1}^q,L)=E_n(C_{\beta,\infty}^{q})_C=
\end{equation*}
\begin{equation}\label{Th_3}
=E_n(C_{\beta,1}^{q})_L=q^n\left(\dfrac{4}{\pi}+\gamma_n\frac{q^{2n}}{1-q^{2n}}\right),
\end{equation}
\emph{де $|\gamma_n|\leqslant\dfrac{16}{3\pi}$.}

\textbf{\emph{Доведення.}} Знайдемо двосторонні оцінки правої частини формули \eqref{dn}. Оскільки,
\begin{equation*}
\left|\sum\limits_{\nu=1}^{\infty}\dfrac{q^{(2\nu+1)n}}{2\nu+1}\sin\left((2\nu+1)\theta_n\pi-\dfrac{\beta\pi}{2}\right)\right|\leqslant
\end{equation*}
\begin{equation*}
\leqslant 
\sum\limits_{\nu=1}^{\infty}\dfrac{q^{(2\nu+1) n}}{2\nu+1}
\leqslant\frac{1}{3}\frac{q^{3n}}{1-q^{2n}},\; n\in\mathbb{N},
\end{equation*}
і в силу формули (64) роботи \cite{My_arxiv_2012} 
\begin{equation*}
1-|\sin(\theta_n\pi-\dfrac{\beta\pi}{2})|\leqslant\frac{q^{2n}}{1-q^{2n}},\; n\in\mathbb{N},
\end{equation*}
то одержуємо для довільних $n\in\mathbb{N}$, $q\in(0,1)$ і $\beta\in\mathbb{R}$
\begin{equation*}
\left|\sum\limits_{\nu=0}^{\infty}\dfrac{q^{(2\nu+1)n}}{2\nu+1}\sin\left((2\nu+1)\theta_n\pi-\dfrac{\beta\pi}{2}\right)\right|\geqslant
1-\left(1-|\sin(\theta_n\pi-\dfrac{\beta\pi}{2})|\right)
-
\end{equation*}
\begin{equation*}
-\left|\sum\limits_{\nu=1}^{\infty}\dfrac{q^{(2\nu+1)n}}{2\nu+1}\sin\left((2\nu+1)\theta_n\pi-\dfrac{\beta\pi}{2}\right)\right|
\geqslant
\end{equation*}
\begin{equation}\label{ots1}
\geqslant
q^n\left(1-\frac{4}{3}\frac{q^{2n}}{1-q^{2n}}\right),
\end{equation}
\begin{equation*}
\left|\sum\limits_{\nu=0}^{\infty}\dfrac{q^{(2\nu+1)n}}{2\nu+1}\sin\left((2\nu+1)\theta_n\pi-\dfrac{\beta\pi}{2}\right)\right|\leqslant
1+\left(1-|\sin(\theta_n\pi-\dfrac{\beta\pi}{2})|\right)
+
\end{equation*}
\begin{equation*}
+\left|\sum\limits_{\nu=1}^{\infty}\dfrac{q^{(2\nu+1)n}}{2\nu+1}\sin\left((2\nu+1)\theta_n\pi-\dfrac{\beta\pi}{2}\right)\right|
\leqslant
\end{equation*}
\begin{equation}\label{ots2}
\leqslant
q^n\left(1+\frac{4}{3}\frac{q^{2n}}{1-q^{2n}}
\right).
\end{equation}

З теореми~2 та оцінок \eqref{ots1} і \eqref{ots2} випливає, що при $n\geqslant n_{q,\beta}$ виконується \eqref{Th_3}. Теорему доведено.

\renewcommand{\refname}{}
\makeatletter\renewcommand{\@biblabel}[1]{#1.}\makeatother

%********************************************************************

\end{document}